\documentstyle[11pt]{article}
\renewcommand{\baselinestretch}{1.5}
\newtheorem{thm}{Theorem}
\newtheorem{cla}{Claim}
\newtheorem{lem}{Lemma}

\newtheorem{cor}{Corollary}

\newtheorem{pro}{Proposition}

\newcommand\proof{\noindent {\bf Proof:\hspace{.05in}}}
\newcommand\qed{\hfill $\Box$ \smallskip}
\newcommand\Qed{\hspace{1cm}$\Box \Box$ \smallskip}

\begin{document}

\title{\bf Rainbow numbers for graphs with cyclomatic number at most two}
\renewcommand\baselinestretch{1.0}
\author{
{\sl Ingo Schiermeyer}
\thanks{Part of this research was performed while the
author was visiting UPJ\v{S} within the project "Research and
Education at UPJ\v{S} - aiming towards excellent  European
universities (EXPERT) under the contract number ITM 26110230056" }
\\
\small Institut f\"ur Diskrete Mathematik und Algebra\\
\small Technische Universit\"at Bergakademie Freiberg
\\ \small09596 Freiberg,
Germany\\
\small Ingo.Schiermeyer@tu-freiberg.de\\
\and {\sl Roman Sot\'ak}\thanks{This work was supported by the
Slovak Science and Technology Assistance Agency under the contract
APVV-0023-10 and Slovak VEGA grant 1/0652/12.
} \\ \small Institute of Mathematics\\
\small P.J. \v{S}af\'arik University in Ko\v{s}ice
\\ \small 04001 Ko\v{s}ice, Slovakia\\
\small Roman.Sotak@upjs.sk}

\date{October 24, 2012}

\maketitle

\renewcommand\baselinestretch{1.5}\normalsize

\maketitle
\renewcommand\baselinestretch{1.2}\normalsize


\begin{abstract}
For a given graph $H$ and $n \geq 1,$ let $f(n, H)$ denote the
maximum number $m$ for which it is possible to colour the edges of
the complete graph $K_n$ with $m$ colours in such a way that each
subgraph $H$ in $K_n$ has at least two edges of the same colour.
Equivalently, any edge-colouring of $K_n$ with at least
$rb(n,H)=f(n,H) + 1$ colours contains a rainbow copy of $H.$ The
numbers $f(n,H)$ and $rb(K_n,H)$ are called {\it anti-ramsey
numbers} and {\it rainbow numbers}, respectively.

In this paper we will classify the rainbow number for a given graph
$H$ with respect to its cyclomatic number. Let $H$ be a graph of
order $p \geq 4$ and cyclomatic number $v(H) \geq 2.$ Then $rb(K_n,
H)$ cannot be bounded from above by a function which is linear in
$n.$ If $H$ has cyclomatic number $v(H) = 1,$ then $rb(K_n, H)$ is
linear in $n.$

We will compute all rainbow numbers for the bull $B,$ which is the
unique graph with $5$ vertices and degree sequence $(1,1,2,3,3).$
Furthermore, we will compute some rainbow numbers for the diamond $D
= K_4-e,$ for $K_{2,3},$ and for the house $H = \overline{P_5}.$

\end{abstract}



\section{Introduction}

We use \cite{BM} for terminology and notation not defined here and
consider finite and simple graphs only. If $K_n$ is edge-coloured in
a given way and a subgraph $H$ contains no two edges of the same
colour, then $H$ will be called a totally multicoloured (TMC) or
rainbow subgraph of $K_n$ and we shall say that $K_n$ contains a TMC
or rainbow $H.$ For a graph $H$ and an integer $n,$ let $f(n,H)$
denote the maximum number of colours in an edge-colouring of $K_n$
with no TMC $H.$ The numbers $f(n,H)$ are called {\it anti-ramsey
numbers} and have been introduced by Erd\H{o}s, Simonovits and S\'os
\cite{ESS}.

We now define $rb(n,H)$ as the minimum number of colours such that
any edge-colouring of $K_n$ with at least $rb(n,H)=f(n,H) + 1$
colours contains a TMC or rainbow subgraph isomorphic to $H.$ The
numbers $rb(n,H)$ will be called {\it rainbow numbers}.


For a given family $\cal H$ of finite graphs
$ext(n,{\cal H}) =:max\{|E(G)| \ | \ H \not\subset G \ \mbox{if} \ H
\in {\cal H}, |V(G)|=n\},$
that is, let $ext(n, {\cal H})$ be the maximum number of edges a
graph $G$ of order $n$ can have if it has no subgraph from $\cal H.$
The graphs attaining the maximum for a given $n$ are called extremal
graphs. The numbers $ext(n, {\cal H})$ are called {\it Tur\'an
numbers} \cite{T}.

For a given graph $H,$ let $\cal H$ be the family of all graphs
which are obtained by deleting one edge from $H.$ If $G$ is a graph
of order $n$ which does not contain any member of $\cal H$ as a
subgraph, then a TMC copy of $G$ and one extra colour for all
remaining edges (of $K_n$) has no TMC subgraph $H.$ Hence, $f(n,H)
\geq ext(n,{\cal H}) + 1.$ Moreover, if we take a rainbow subgraph
with $ext(n, H) + 1$ edges, then it contains a rainbow subgraph
isomorphic to $H.$ Hence (cf. \cite{ESS})

\begin{eqnarray}
ext(n,{\cal H}) + 2 \leq f(n,H) + 1 = rb(n,H) \leq ext(n,H)+1.
\end{eqnarray}

The lower bound is sharp for some graph classes. This has been shown
if $H$ is a complete graph on $k \geq 3$ vertices in \cite{MBNL02,S}
and if $H$ is a matching with $k$ edges and $n \geq 2k+1$ in
\cite{FKSS,S}.

Erd\H{o}s, Simonovits and S\'os \cite{ESS} showed that $f(n,H)/{n
\choose 2} \rightarrow 1 - \frac{1}{d}$ as $n \rightarrow \infty,$
where $d + 1 = min\{\chi(H-e) \ | \ e \in E(H)\},$ and that
$f(n,H) - ext(n, {\cal{H}}) = o(n^2).$ Hence the rainbow numbers
are asymptotically known if $min\{\chi(H-e) \ | \ e \in E(H)\}
\geq 3.$ If $min\{\chi(H-e) \ | \ e \in E(H)\} \leq 2,$ then the
situation is quite different.

For cycles the following result (which has been conjectured by Erd\H{o}s, Simonovits and S\'os \cite{ESS})
has been shown by Montellano-Ballesteros and Neumann-Lara \cite{MBNL05}.

\begin{thm}\label{t1} \cite{MBNL05} \
Let $n \geq k \geq 3.$ Then
$rb(n,C_k)=\lfloor\frac{n}{k-1}\rfloor{k-1 \choose 2} + {r \choose
2} + \lceil\frac{n}{k-1}\rceil,$ where $r$ is the residue of $n$
modulo $k-1.$
\end{thm}

Gorgol \cite{G} has considered a cycle $C_k$ with a pendant edge, denoted $C_k^+,$ and computed all rainbow numbers.

\begin{thm}
\cite{G} \\
$rb(n,C_k^+) = rb(n,C_k),$ for $n \geq k+1.$
\end{thm}

However, if we add two (or more) edges to a cycle $C_k,$ the situation becomes
surprisingly interesting.

\section{Rainbow numbers and cyclomatic number}

Before presenting our main results we introduce some additional
notation (similar as in \cite{MB}), which will be used frequently in
the proofs.

For a set $C$ of colours let $\Gamma: E(G) \rightarrow C(G)$ be an
edge-colouring of $G$ and let $F$ be a subgraph of $G.$ $\Gamma[F]$
will denote the image of $\Gamma$ restricted to the set $E(F),$ and
$F$ will be called {\it rainbow}, if no pair of edges of $F$ receive
the same colour., i.e. $|E(F)| = |\Gamma[F]|.$ Given $w \in V(F),$
let $S(w,F)$ be the set of colours $c \in \Gamma(F)$ such that for
every $e \in E(F),$ if $\Gamma(E) = c$ then $e$ is incident to the
vertex $w.$ If $F \cong G,$ then we will write $S(w).$ Observe that
for every $w \in V(F),$ if $c \in S(w,F)$ then

\begin{eqnarray}\label{f0}
c \notin \Gamma[F \setminus \{w\}] \ \mbox{and} \ |\Gamma[F
\setminus \{w\}]| = |\Gamma[F]| - |S(w,F)|.
\end{eqnarray}

For each colour $i$ let $E_i = \{e \in E(G) | c(e) = i\}.$ If $i \in
S(v)$ for a vertex $v \in V(G),$ then $i$ is called a {\it unique
colour at vertex v}. Let $C^* = \cup_{v \in V(G)}S(v)$ be the set of
all unique colours. So $G[E_i] \cong K_{1,|E_i|}$ for each $i \in
C^*.$

The {\it cyclomatic number} $v(G)$ of a connected graph $G$ is
defined as $v(G) = |E(G)| - |V(G)| + 1$ and measures, how many edges
have to be deleted from a given graph $G$ in order to make it
acyclic.


\begin{thm}
Let $H$ be a graph of order $p \geq 4$ and cyclomatic number $v(H)
\geq 2.$ Then $rb(K_n, H)$ cannot be bounded from above by a
function which is linear in $n.$
\end{thm}

For the proof of this Theorem we apply the following Theorem of
Erd\H{o}s.

\begin{thm}(Erd\H{o}s \cite{E})\label{tE}\\
For every pair of integers $g \geq 4, k \geq 3$ there is a graph
with girth at least $g$ and with chromatic number at least $k.$
\end{thm}

\proof Choose $g = p+1$ and $k$ arbitrarily. Then by Theorem
\ref{tE} there is a graph $F$ with girth $g \geq p+1$ and $\chi(F)
\geq k+1.$ We may assume $\delta(F) \geq k,$ since otherwise we
could delete vertices of degree $\leq k-1.$ For $t \geq 1$ let $F_t$
be the graph consisting of $t$ pairwise vertex-disjoint copies of
$F.$ Let $n = t|V(F)|$ and choose a copy of $F_t \subset K_n.$
Colour all edges of $F_t$ distinct and all remaining edges with
another colour $c^*.$ Then $|C(K_n)| \geq \frac{kn}{2} + 1 >
\frac{k}{2} \cdot n.$ Since $F_t$ has girth $g > p$ and all cycles
in $H$ have length at most $p,$ we have to delete at least two edges
from $H$ in order to make it acyclic. Hence in every copy of $H$ at
least two edges have the same colour $c^*$  and so $rb(K_n, H) \geq
\frac{kn}{2} + 2.$ \qed

If $v(H)=1,$ then the situation is quite different. If $H$ has
cyclomatic number $v(H)=1,$ then we can show an upper bound for
$rb(K_n,H),$ which is linear in $n.$ A lower bound for those rainbow
number depends on the length of the unique cycle of $H.$ Therefore
we will prove the following Theorem.

\begin{thm}\label{t5}
Let $H$ be a graph of order $p \geq 5$ and cyclomatic number $v(H) =
1.$ If $H$ contains a cycle with $k$ vertices, $3 \leq k \leq p-2,$
then
$$\lfloor\frac{n}{k-1}\rfloor{k-1 \choose 2}
+ {r \choose 2} + \lceil\frac{n}{k-1}\rceil \leq rb(K_n, H) \leq
(p-2)n - p \cdot \frac{p-3}{2},$$ where $n \geq p$ and $r$ is the
residue of $n$ modulo $k-1.$
\end{thm}

\proof The lower bound is given by Theorem \ref{t1}. Our proof for
the upper bound is by induction on the order $n.$ Let $G \cong K_n$
with $|C(G)| = (p-2)n - p \cdot \frac{p-3}{2}.$

If $n = p,$ then $(p-2)p - p \cdot \frac{p-3}{2} = {p \choose 2}.$
Clearly, if all ${p \choose 2}$ edges are coloured distinct, there
is a rainbow copy of $H.$ So we may assume that $n \geq p+1.$ If
there is a vertex $v$ with $|S(v)| \leq p-2,$ then we apply
induction. Hence we may assume that $|S(v)| \geq p-1$ for every
vertex $v \in V(G).$ By Theorem \ref{t1} $G$ contains a rainbow
cycle $C_k.$ Now we can extend this cycle $C_k$ to a copy of $H$ by
properly choosing $p-k$ additional vertices, since $|S(v)| \geq p-1$
for every vertex $v \in V(G).$ \qed

\section{Rainbow numbers for the bull}

Let $B$ be the unique graph with $5$ vertices and degree sequence
$(1,1,2,3,3),$ which is called the {\it bull}. By Theorem \ref{t5}
we obtain $n \leq rb(K_n,B) \leq 3n-5.$ Here we have been able to
compute all rainbow numbers for the bull.

\begin{thm}
$rb(5,B)=6$ and $rb(n,B) = n+2$ for $n \geq 6.$
\end{thm}

\begin{pro}
$rb(K_5,B) = 6.$
\end{pro}

\proof Let $G \cong K_5$ and choose a $K_3$ and a $K_2,$ which are
vertex disjoint. Colour the four included edges with four distinct
colours and all remaining edges with a fifth colour. Then $G$
contains no rainbow bull. So consider an edge colouring of $G$ with
six colours. With $rb(K_5,K_3)=5$ we conclude that $G$ contains a
rainbow $K_3.$ We distinguish three cases.

\noindent{\bf Case 1} \ $G$ contains a rainbow $K_4$\\
Let $w_1, w_2, w_3, w_4$ be the vertices of the $K_4,$ whose six
edges are coloured with colours $1, 2, \ldots, 6.$ Let $w_5$ be the
fifth vertex and consider an edge incident with $w_5,$ say $w_5w_1.$
Then $c(w_1w_5) \in \{1, 2, \ldots, 6\}$ and we always find a
rainbow bull.

\medskip

\noindent{\bf Case 2} \ $G$ contains a rainbow $K_4-e,$ but no rainbow $K_4$\\
Let $w_1, w_2, w_3, w_4$ be the vertices of the $K_4-e,$ whose five
edges are coloured with colours $1, 2, \ldots, 5.$ Then there is an
edge $w_5w_i$ for some $1 \leq i \leq 4$ with $c(w_5w_i) = 6$ and we
always find a rainbow bull.

\medskip

\noindent{\bf Case 3} \ $G$ contains a rainbow $K_3,$ but no rainbow $K_4-e$\\
Let $w_1, w_2, w_3$ be the vertices of the rainbow $K_3,$ whose
edges are coloured with the three colours $1, 2, 3,$ and let $w_4,
w_5$ be the two other vertices. If (i.e.) $c(w_1w_4) = 4$ and
$c(w_2w_5) = 5,$ then there is rainbow bull. Repeating these
arguments we either find a rainbow bull or conclude that
$c(w_4w_1)=4, c(w_5w_1)=5, c(w_4w_5)=6.$ Now any colour from $\{1,
2, \ldots, 6\}$ for an edge $w_iw_j$ for $2 \leq i \leq 3$ and $4
\leq j \leq 5$ yields a rainbow bull. \qed

\begin{lem}
$rb(K_n,B) \geq n+2$ for all $n \geq 6.$
\end{lem}

\proof Let $r$ be the remainder of $n$ modulo $3.$ Then $n =
3(\lfloor\frac{n}{3}\rfloor - r) + 4r.$ Now choose a collection of
pairwise vertex disjoint $(\lfloor\frac{n}{3}\rfloor - r) \ C_3$'s
and $r \ C_4$'s and colour all their edges with $n$ distinct colours
and all remaining edges of $K_n$ with another colour. Then $G$
contains no rainbow bull. \qed

\begin{pro}
$rb(K_6,B) = 8.$
\end{pro}

\proof Let the edges of $K_6$ be coloured with $8$ colours. If
$|S(v)| \leq 2$ for a vertex $v \in V(G),$ then there is a rainbow
bull in $G - v$ by induction. Hence we may assume $|S(v)| \geq 3$
for every vertex $v \in V(G).$ Then $|C| \geq \lceil\frac{6 \cdot
3}{2}\rceil = 9,$ a contradiction. \qed

\begin{pro}
$rb(K_n,B) = rb(K_{n-1},B) +1$ for all $n \geq 7.$
\end{pro}

\proof Let $G \cong K_n$ with $|C(G)| = rb(K_{n-1},B) +1.$ If
$|S(v)| \leq 1$ for a vertex $v \in V(G),$ then there is a rainbow
bull in $G - v$ by induction. Hence we may assume

\begin{eqnarray}\label{f1}
|S(v)| \geq 2 \ \mbox{for every} \ v \in V(G).
\end{eqnarray}

\begin{lem}\label{l2}
$|S(v)| = 2$ for every vertex $v \in V(G).$
\end{lem}

\proof Suppose $|S(w)| \geq 3$ for a vertex $w \in V(G).$ Let
$c(ww_i) = i \in C^*$ for three vertices $w_1, w_2, w_3.$ If
$c(w_iw_j) = p$ for some $i,j \in \{1,2,3\}$ with $i \neq j$ and $p
\notin \{1,2,3\},$ then $c(w_iu) = c(w_ju) = p$ for all $u \in V
\setminus \{w, w_1, w_2, w_3\},$ otherwise we would have a rainbow
bull generated by $w, w_1, w_2, w_3, u.$ As $V \setminus \{w, w_1,
w_2, w_3\} \neq \emptyset,$ we conclude that
$c(w_1w_2)=c(w_1u)=c(w_2u)=c(w_2w_3)=c(w_1w_3)=p.$ Then $S(w_1)
\subseteq \{1\},$ contradicting (\ref{f1}). \qed

Let $c_1 = |\{i \in C^* | |E_i| = 1\}|, c_2 = |\{i \in C^* | |E_i|
\geq 2\}|$ and $c_3 = |C| - (c_1 + c_2).$ Suppose $|C| = n+2.$ Then
by Lemma \ref{l2}

\begin{eqnarray}\label{f2}
2c_1 + c_2 & = & 2n\\
c_1 + c_2 + c_3 & = & n+2
\end{eqnarray}

\begin{lem}\label{l3}
Let $H \subset G$ be a subgraph with $V(H) = \{w_1, w_2, w_3\}$ and
$c(w_iw_{i+1})=i$ (indices taken mod 3) for $i = 1,2,3.$ If $|E_i| =
1$ for $i = 1,2,3,$ and $G$ contains no rainbow bull, then there is
some $k \in C \setminus \{1,2,3\}$ such that $c(e) = k$ for all
edges $e = uw$ with $u \in V(H), w \in V(G) \setminus V(H).$
\end{lem}

\proof Choose a vertex $u \in V(G) \setminus V(H).$ Then $c(w_1u) =
k$ for some $k \notin \{1,2,3\},$ say $k = 4.$ Since there is no
rainbow bull we conclude first that $c(vw_i) = 4$ for all vertices
$v \in V(G) \setminus (V(H) \cup \{u\})$ and $i = 2,3,$ then
$c(uw_i) = 4$ for $i = 2,3,$ and finally $c(vw_1) = 4$ for all
vertices $v \in V(G) \setminus (V(H) \cup \{u\}).$ \qed

\begin{lem}\label{l4}
Let $H \subset G$ be a subgraph with $V(H) = \{w_1, w_2, w_3, w_4\}$
and $c(w_iw_{i+1})=i$ (indices taken mod 4) for $i = 1,2,3,4.$ If
$|E_i| = 1$ for $i = 1,2,3,4,$ and $G$ contains no rainbow bull,
then there is some $k \in C \setminus \{1,2,3,4\}$ such that $c(e) =
k$ for all edges $e = uw$ with $u \in V(H), w \in V(G) \setminus
V(H),$ and $c(w_1w_3) = c(w_2w_4) = k.$
\end{lem}

\proof Choose a vertex $u \in V(G) \setminus V(H).$ Then $c(w_1u) =
k$ for some $k \notin \{1,2,3,4\},$ say $k = 5.$ Since there is no
rainbow bull we conclude first that $c(w_2w_4) = 5,$ then $c(vw_i) =
5$ for all vertices $v \in V(G) \setminus V(H)$ and $i =1, 2, 3, 4,$
and finally $c(w_1w_3) = 5.$ \qed

\begin{lem}\label{l5}
For $k \geq 5$ let $H=P_k \subset G$ be a subgraph with $V(H) = \{w_1, w_2, \ldots, w_k\}$ and
$c(w_iw_{i+1})=i$ for $1 \leq i \leq k-1.$ If $|E_i| = 1$ for $1 \leq i \leq k-1,$
then $G$ contains a rainbow bull.
\end{lem}

\proof We consider the edge $w_2w_4.$ By the assumption $c(w_2w_4)
\neq c(w_iw_{i+1})$ for $1 \leq i \leq k-1.$ Hence $\{w_1, w_2, w_3,
w_4, w_5\}$ generate a rainbow bull. \qed

\medskip

\noindent{\bf Case 1} $c_1 = n, c_2 = 0, c_3 = 2$\\
Let $H$ be the graph induced by all colours $i$ with $|E_i|=1.$ So
$|E(H)| = n.$ Since $|S(v)| = 2$ for every vertex $v \in V(G)$ and
by Lemma \ref{l5} we conclude that $H$ is a collection of pairwise
vertex-disjoint $3$-cycles and $4$-cycles. Moreover, $n \geq 7,$ so
with Lemma \ref{l3} and Lemma \ref{l4} we conclude $|C| = n+1 <
n+2,$ a contradiction.

\medskip

\noindent{\bf Case 2} $c_1 = n-1, c_2 = 2, c_3 = 1$\\
Let $H$ be the graph induced by all colours $i$ with $|E_i|=1.$ So
$|E(H)| = n-1.$ Since $|S(v)| = 2$ for every vertex $v \in V(G)$ and
by Lemma \ref{l5} we conclude that $H$ is a collection of pairwise
vertex-disjoint $3$-cycles and $4$-cycles and a path $P_k$ for some
$1 \leq k \leq 4.$ For $k=1$ using Lemma \ref{l3} and Lemma \ref{l4}
we obtain $|C| = n-1 + 1 = n < n+2,$ a contradiction. Hence we may
assume that $2 \leq k \leq 4.$ Let $w_1, \ldots, w_k$ be the
vertices of the path with edges $w_iw_{i+1}$ for $1 \leq i \leq
k-1.$ Suppose $G$ contains no rainbow bull. By Lemma \ref{l3} and
Lemma \ref{l4} there is a colour $p$ such that $c(w_iv)=p$ for all
edges $w_iv$ with $1 \leq i \leq k$ and $v \in V(G) \setminus
V(P_k).$ Hence $p \notin S(w_1) \cup S(w_k).$ With $|S(w_1)| =
|S(w_k)| = 2$ and $c(w_1w_k) \notin C_1$ we conclude that $k = 4.$
Furthermore, using $c_2 = 2,$ we obtain $c(w_1w_3) = c(w_1w_4)$ and
$c(w_1w_4) = c(w_2w_4),$ a contradiction. \qed

\medskip

\noindent{\bf Case 3} $c_1 = n-2, c_2 = 4, c_3 = 0$\\
Let $|E_i| \geq 2$ for $i = 1,2,3,4$ and $i \in C^*.$

\begin{cla}
$\sum_{i=1}^4|E(i)| \leq 4n-12$
\end{cla}

\proof First suppose that there are two vertices $w_1, w_2 \in V(G)$
such that $1,2 \in S(w_1)$ and $3,4 \in S(w_2).$ Then
$\sum_{i=1}^4|E(i)| \leq 2(n-2) + 1 = 2n-3 < 4n-12$ for $n \geq 7.$

Suppose next that there are three vertices $w_1, w_2, w_3 \in V(G)$
such that $1,2 \in S(w_1), 3 \in S(w_2), 4 \in S(w_3).$ Now $w_2$ is
incident with some edge $e$ with $c(e) \notin \{1,2,3,4\},$ thus
$\sum_{i=1}^4|E(i)|$ $ \leq 3(n-3) + 3-1 = 3n-7 < 4n-12$ for $n \geq
7.$

Finally suppose that there are four vertices $w_1, w_2, w_3, w_4 \in
V(G)$ such that $i \in S(w_i)$ for $1 \leq i \leq 4.$ Now each $w_i$
is incident with some edge $e_i$ with $c(e_i) \notin \{1,2,3,4\},$
thus $\sum_{i=1}^4|E(i)| \leq 4(n-4) + 6- \frac{4}{2} = 4n-12.$ \qed

Now we obtain
\begin{eqnarray}
|E(G)| = {n \choose 2} & \leq & (n-2) + 4n-12\\
\Leftrightarrow n^2 - 11n + 28 & \leq & 0\\
\Rightarrow 4 \leq n & \leq & 7
\end{eqnarray}

By the assumption $(n \geq 7)$ we consider the case $n=7.$ Moreover,
if $\sum_{i=1}^4|E(i)| < 4n-12,$ then $n < 7,$ a contradiction.

Thus $\sum_{i=1}^4|E(i)| = 4n-12$ and $c_1 = 5$ and $c_2 = 4.$
Following the previous arguments there are vertices $w_1, w_2,
\ldots, w_7$ with $w_1w_2, w_3w_4, w_5w_6, w_6w_7, w_5w_7 \in E(G)$
and $c(w_1w_2)=5, c(w_3w_4)=6, c(w_5w_6)=7, c(w_6w_7)=8,
c(w_5w_7)=9$ such that $i \in C^*$ for $5 \leq i \leq 9.$ By Lemma
\ref{l3} we conclude that there is only one colour among all edges
$w_iw_j$ for $1 \leq i \leq 4$ and $5 \leq j \leq 7.$ But we have
$c_3 = 0,$ a contradiction. \qed \Qed

\section{Rainbow numbers for the diamond}

We now consider the graph $D = K_4-e,$ which is called the {\it
diamond}. This graph contains a $C_3$ and has cyclomatic number
$v(D)=2.$

In \cite{GKL} projective planes have been used to construct an
infinite family of graphs with girth 6 (i.e. $\{C_3, C_4,
C_5\}$-free) having $\Omega(n^{\frac{3}{2}})$ edges. This means that
$ext(n, \{C_3, C_4, C_5\}) = \Omega(n^{\frac{3}{2}})$ and so by
(\ref{f0}) we deduce

\begin{cor}\label{c1}
$rb(K_n, D) = \Omega(n^{\frac{3}{2}}).$
\end{cor}

Montellano-Ballesteros \cite{MB} has shown an upper bound for the
rainbow number of the diamond.

\begin{thm}\label{t7}
\cite{MB} \ For every $n \geq 4,$
\begin{eqnarray}\label{f7}
ext(n,\{C_3,C_4\})+2  \leq rb(K_n,D) \leq ext(n,\{C_3,C_4\})+(n+1).
\end{eqnarray}
\end{thm}

For the extremal numbers $ext(n,\{C_3,C_4\})$ the following bounds
have been shown in \cite{GKL}.

\begin{thm}\label{t8}
$ext(n,\{C_3,C_4\}) \leq \frac{1}{2}n\sqrt{n-1}.$
\end{thm}

\begin{thm}\label{t9}
$$\frac{1}{2\sqrt{2}} \leq lim \ inf_{n \rightarrow \infty}
\frac{ext(n,\{C_3,C_4\})}{n^{\frac{3}{2}}} \leq lim \ sup_{n
\rightarrow \infty} \frac{ext(n,\{C_3,C_4\})}{n^{\frac{3}{2}}} \leq
\frac{1}{2}.$$
\end{thm}

By (\ref{f7}) we obtain

\begin{cor}\label{c2}
$rb(K_n,D) \leq \frac{1}{2}n\sqrt{n-1} + (n+1).$
\end{cor}

Corollary \ref{c1} and Corollary \ref{c2} give

\begin{cor}\label{c3}
$rb(K_n,D) = \Theta(n^{\frac{3}{2}}).$
\end{cor}

With Theorem \ref{t7} and Theorem \ref{t9} we obtain

\begin{cor}
$$\frac{1}{2\sqrt{2}} \leq lim \ inf_{n \rightarrow \infty}
\frac{rb(K_n,D)}{n^{\frac{3}{2}}} \leq lim \ sup_{n \rightarrow
\infty} \frac{rb(K_n,D)}{n^{\frac{3}{2}}} \leq \frac{1}{2}.$$
\end{cor}

In the following we will compute the rainbow numbers $rb(K_n,D)$ for
$4 \leq n \leq 10.$ It turns out that these values are all equal to
the lower bound in (\ref{f7}). The extremal numbers
$ext(n,\{C_3,C_4\}$ have been computed in \cite{GKL}. Here we list
the values for $4 \leq n \leq 16.$

\medskip

\begin{center}
\begin{tabular}{c|c|c|c|c|c|c|c|c|c|c|c|c|c}
n & 4 & 5 & 6 & 7 & 8 & 9 & 10 & 11 & 12 & 13 & 14 & 15 & 16\\
\hline
$ext(n,\{C_3,C_4\}$ & 3 & 5 & 6 & 8 & 10 & 12 & 15 & 16 & 18 & 21 & 23 & 26 & 28 \\
\end{tabular}
\end{center}

\begin{thm}\label{t10}
$rb(K_n,D) = ext(n,\{C_3,C_4\}) +2$ for $4 \leq n \leq 10.$
\end{thm}

\begin{cor}
$ext(n,\{K_{1,3}+e,C_4\}) = ext(n,\{C_3,C_4\})$ for $4 \leq n \leq
10.$
\end{cor}

\proof ${\cal H(D)} = \{K_{1,3}+e,C_4\}$ and so $
ext(n,\{C_3,C_4\})+2 \leq ext(n,\{K_{1,3}+e,C_4\})+2 \leq rb(K_n,D)
= ext(n,\{C_3,C_4\})+2.$ \qed

In the proof of Theorem \ref{t10} we will use the following helpful
Lemma.

\begin{lem}
$rb(K_n, D) \leq rb(K_n, K_{2,3})$ for all $n \geq 5.$
\end{lem}

\proof Let $K_n$ be edge coloured and let $F$ be a rainbow $K_{2,3}$
with vertices $w_1, w_2, \ldots, w_5$ and edges $w_iw_j$ for $1 \leq
i \leq 2, 3 \leq j \leq 5.$ Let $c(w_1w_{i+2})=i, c(w_2w_{i+2})=i+3$
for $1 \leq i \leq 3.$ Now consider the edge $w_1w_2.$ If $c(w_1w_2)
= x$ with $x \in C \setminus \{1,2, \ldots, 6\},$ then $F + w_1w_2$
contains a rainbow $D.$ Otherwise, if $c(w_1w_2) = c(e) \in \{1,2,
\ldots,6\}$ for an edge $e \in E(F),$ then $F + w_1w_2 - e$ contains
a rainbow $D.$ \qed

\medskip

\noindent\proof (of Theorem \ref{t10})\\
We will give a seperate proof for each value of $n, 4 \leq n \leq
10.$\\
\noindent{$\bf rb(K_4,D) = 5$}\\
\proof Let the edges of $K_4$ be coloured with $5$ colours. Choose $5$ edges with
distinct colours in $K_4.$ \qed

\medskip

\noindent{$\bf rb(K_5,D) = 7$}\\
\proof Let the edges of $K_5$ be coloured with $7$ colours.
If $|S(v)| \leq 2$ for a vertex $v \in V(G),$ then there is a rainbow diamond in $G - v$ by induction.
Hence we may assume $|S(v)| \geq 3$ for every vertex $v \in V(G).$ Then $|C| \geq \lceil\frac{5 \cdot 3}{2}\rceil = 8,$
a contradiction. \qed

\medskip

\noindent{$\bf rb(K_6,D) = 8$}\\
\proof Let the edges of $K_6$ be coloured with $8$ colours and let
$F$ be a rainbow subgraph of $K_6$ with $8$ edges. Since $rb(K_6,
K_{1,4}+e)=8,$ (cf. \cite{GL}) we may assume that $F$ contains a
$K_{1,4}.$ Suppose that $F$ contains a $K_{1,5}$ with vertices $w,
w_1, w_2, w_3, w_4, w_5$ and edges $ww_i$ for $1 \leq i \leq 5.$
Then $F[N(w)]$ contains three edges and hence $F$ contains a
diamond. So there is a rainbow diamond. So we may assume that $F$
contains a $K_{1,4}$ (and no $K_{1,5}$) with vertices $w, w_1, w_2,
w_3, w_4$ and edges $ww_i$ for $1 \leq i \leq 4.$ Let $w_5$ be the
sixth vertex with $w_5 \notin N_F(w).$ If $F[N(w)]$ contains at
least three edges, then there is a rainbow diamond. If $d_F(w_5)
\geq 3,$ then there is a rainbow $K_{2,3}$ and thus a rainbow
diamond. So we may assume that $F[N(w)]$ contains exactly two edges,
say $w_1w_2, w_3w_4,$ and that $d_F(w_5)=2.$ If $w_5w_i, w_5w_{i+1}
\in E(F)$ for $i = 1$ or $i = 3,$ then there is a rainbow diamond.
So we may assume that $w_5w_1, w_5w_3 \in E(F).$ Now consider the
edge $ww_5.$ Then either $\{w,w_1,w_2,w_5\}$ or $\{w,w_3,w_4,w_5\}$
generate a rainbow diamond. \qed

\medskip

\noindent{$\bf rb(K_7,D) = 10$}\\
\proof Let the edges of $K_7$ be coloured with $10$ colours. If $|S(v)| \geq 3$ for every vertex
$v \in V(G),$ then $|C| \geq \lceil\frac{7 \cdot 3}{2}\rceil = 11,$ a contradiction. Hence there is a vertex
$v \in V(G)$ with $|S(v)| \leq 2.$ Then $G - v$ contains a rainbow diamond by induction. \qed

\medskip

\noindent{$\bf rb(K_8,D) = 12$}\\
\proof Let the edges of $K_8$ be coloured with $12$ colours. If
$|S(v)| \leq 2$ for a vertex $v \in V(G),$ then there is a rainbow
diamond in $G - v$ by induction. Hence we may assume $|S(v)| \geq 3$
for every vertex $v \in V(G).$ With $12 = |C| \geq \lceil\frac{8
\cdot 3}{2}\rceil = 12$ we deduce $|S(v)| = 3$ for every vertex $v
\in V(G).$ Then there are four vertices $w,w_1,w_2,w_3 \in V(G)$
such that $c(ww_i)=i \in S(w)$ and $i \in C^*$ for $1 \leq i \leq
3.$ If $c(w_1w_2) \neq c(w_1w_3)$ or $c(w_1w_2) \neq c(w_2w_3),$
then there is rainbow diamond. Hence we may assume that $c(w_1w_2) =
c(w_1w_3) = c(w_2w_3) = p,$ say $p=4,$ and therefore $p \notin
S(w_i).$ Since $|S(w_i)| = 3$ for $1 \leq i \leq 3$ we conclude that
there are vertices $w_4, w_5, \ldots, w_9$ such that $w_1w_4,
w_1w_5, w_2w_6, w_2w_7, w_3w_8, w_3w_9 \in E(G)$ and $c(w_1w_i)=i+1$
for $4 \leq i \leq 5,$ $c(w_2w_i)=i+1$ for $6 \leq i \leq 7,$ and
$c(w_3w_i)=i+1$ for $8 \leq i \leq 9.$ Moreover, $\{w_4,w_5\},
\{w_6,w_7\},\{w_8,w_9\}$ are pairwise vertex-disjoint, since
otherwise $w, w_1, w_2, w_3, w_i$ for a vertex $w_i, 4 \leq i \leq
9,$ would create a rainbow diamond. But then $|V(G)| \geq 10,$ a
contradiction. \qed

\medskip

\noindent{$\bf rb(K_9,D) = 14$}\\
\proof Let the edges of $K_9$ be coloured with $14$ colours.
We can follow the proof for $rb(K_8,D) = 12$ in all steps and obtain the same final contradiction. \qed

\medskip

\noindent{$\bf rb(K_{10},D) = 17$}\\
\proof Let the edges of $K_{10}$ be coloured with $17$ colours. If
$|S(v)| \leq 3$ for a vertex $v \in V(G),$ then there is a rainbow
diamond in $G - v$ by induction. Hence we may assume $|S(v)| \geq 4$
for every vertex $v \in V(G).$ Then $|C| \geq \lceil\frac{10 \cdot
4}{2}\rceil = 20,$ a contradiction. \qed \Qed

\section{Rainbow numbers for $K_{2,3}$}

In \cite{AJ} anti-ramsey numbers have been computed for small
complete bipartite graphs.

\begin{thm}\label{t11}
$f(K_n,K_{2,t}) = \frac{\sqrt{t-2}}{2}n^{\frac{3}{2}} +
O(n^{\frac{4}{3}}).$
\end{thm}

So we deduce

\begin{cor}
$rb(K_n, K_{2,3}) = \Theta(n^{\frac{3}{2}}).$
\end{cor}

\medskip
We say that $\{w_1, w_2, w_3, w_4, w_5\}$ generate a $K_{2,3},$
where $d(w_1)=d(w_2)=3$ and $d(w_3)=d(w_4)=d(w_5)=2.$ Now we will
compute the rainbow numbers $rb(K_n, K_{2,3})$ for $5 \leq n \leq
8.$
\medskip

\noindent{$\bf rb(K_5, K_{2,3}) = 8$}\\
\proof Consider a $K_4$ and a $K_1,$ which are vertex disjoint. Colour the six
edges of the $K_4$ distinct and all remaining edges with another colour. Then there is no
rainbow $K_{2,3}.$ Hence $rb(K_5, K_{2,3}) \geq 8.$ Consider $8$ edges of a $K_5$ with distinct
colours. The induced graph is isomorphic to either $K_5 - P_3$ or $K_5 - 2K_2.$
In both cases we find a rainbow $K_{2,3}.$ \qed

\noindent{$\bf rb(K_6, K_{2,3}) = 10$}\\
\proof Consider a $K_4$ with vertices $w_1, w_2, w_3, w_4,$ and let
$w_5, w_6$ be the two other vertices. Colour all edges of the $K_4$
distinct with colours $1,2, \ldots,6$ and let $c(w_5w_6)=7,
c(w_5w_i)=8, c(w_6w_i)=9$ for $1 \leq i \leq 4.$ Then there is no
rainbow $K_{2,3}.$ Hence $rb(K_6, K_{2,3}) \geq 10.$

Let the edges of $K_6$ be coloured with $10$ colours and let $F$ be
a rainbow subgraph of $K_6$ with $10$ edges. Since $rb(K_6,D)=8,$
there is a rainbow $K_4-e.$ If there is a rainbow $K_4,$ then $|N(v)
\cap V(K_4)| \geq 2$ for a vertex $v$ outside of the $K_4$ and we
obtain a rainbow $K_{2,3}.$ Hence we may assume that $F$ contains a
rainbow $K_4-e,$ but no rainbow $K_4.$ Let $w_1,w_2,w_3,w_4$ be the
vertices of the $K_4-e$ with edges
$w_1w_2,w_2w_3,w_3w_4,w_4w_1,w_1w_3.$ If $w_5w_i, w_5w_{i+2} \in
E(F)$ for some $i, 1 \leq i \leq 2,$ then there is a rainbow
$K_{2,3}.$ Hence we may assume that $d_{K_4-e}(w_5) = d_{K_4-e}(w_6)
= 2$ and $w_5w_6 \in E(F).$ Moreover, $N_{K_4-e}(w_5) =
\{w_i,w_{i+1}\}$ and $N_{K_4-e}(w_6) = \{w_j,w_{j+1}\}$ for some $1
\leq i,j \leq 4.$ If $j=i$ or $j = i+2,$ then there is a rainbow
$K_{2,3}.$ So assume $j = i+1.$ If $i=1, j=2$ or $i=3, j=4,$ then
there is a rainbow $K_{2,3}.$ Hence assume $i=2, j=3.$ Then there is
rainbow $W_5$ (wheel on six vertices). We relabel the vertices $w_1,
w_2, \ldots, w_6$ such that $w_iw_{i+1}, w_iw_6 \in E(F)$ for $1
\leq i \leq 5$ ($i+1$ modulo 5).

Let $F' = F + w_2w_5.$ Then the three sets of vertices $V_1 = \{w_5,
w_6, w_1, w_2, w_4\},$ $V_2 = \{w_2, w_6, w_1, w_3, w_5\},$ $V_3 =
\{w_3, w_5, w_2, w_4, w_6\}$ each generate a $K_{2,3}.$ Using $V_1,
V_2, V_3$ we observe that $F' - e = F + w_2w_5 - e$ contains a
rainbow $K_{2,3}$ for any edge $e \in E(F).$ \qed

\noindent{$\bf rb(K_7, K_{2,3}) = 12$}\\
\proof Consider a $K_3$ and a $K_4,$ which are vertex disjoint. Let the vertices
of the $K_3$ be $w_1, w_2, w_3,$ and the vertices of the $K_4$ be
$w_4, w_5, w_6, w_7.$ Colour all edges of the $K_3$ distinct with colours $1, 2, 3$ and
all edges of the $K_4$ distinct with colours $4, 5, \ldots, 9.$
Next let $c(w_1w_i)=10, c(w_2w_i)=c(w_3w_i)=11$ for $4 \leq i \leq 7.$
Then $G$ has no rainbow $K_{2,3}.$ Hence $rb(K_7, K_{2,3}) \geq 12.$

Let the edges of $K_7$ be coloured with $12$ colours and let $F$ be
a rainbow subgraph of $K_7$ with $12$ edges such that $\Delta(F)$ is
maximum. If $|S(v)| \leq 2$ for a vertex $v \in V(G),$ then $G$ has
a rainbow $K_{2,3}$ by induction. Hence we may assume that $|S(v)|
\geq 3$ for every vertex $v \in V(G).$ With $|E(F)| \geq
\lceil\frac{3 \cdot 7}{2}\rceil = 11
> \frac{3 \cdot 7}{2}$ we deduce that $4 \leq \Delta(F) \leq 6.$ We
distinguish these three cases.

\noindent{\bf Case 1} $\Delta(F) = 6$\\
Let $w \in V(F)$ be a vertex with $d_F(w) = 6.$ Since $|S(v)| \geq 3
= 2 + 1$ for all vertices $v \in V(F) \setminus \{w\},$ we conclude
that $|C^*| = c_1 + c_2 \geq 6 + \frac{6 \cdot 2}{2} = 12$ and $2c_1
+ c_2 \geq 6 + 6 \cdot 3 = 24$ implying $c_1 = 12, c_2 = 0.$ But
then $c_ 3 \geq 1,$ since $|E(K_7)| = 21 > 12,$ and so $|C| \geq
13,$ a contradiction. \qed

\noindent{\bf Case 2} $\Delta(F) = 5$\\
Let $w_1 \in V(F)$ be a vertex with $d_F(w_1) = 5$ and let $w_2,
\ldots, w_6$ be its neighbours. Let $w_7$ be the seventh vertex.
With $|S(w)| \geq 3$ we obtain a rainbow $K_{2,3}$ with vertices
$\{w_1, w_7, x,y, z\},$ where $x,y,z \in \{w_2, \ldots, w_6\}.$ \qed

\noindent{\bf Case 3} $\Delta(F) = 4$\\
We first show that there is a vertex $v \in V(G)$ such that $|S(v)|
\geq 4.$

Suppose that $|S(v)|=3$ for all vertices $v \in V(G).$ Then $2c_1 +
c_2 = 21$ and $c_1 + c_2 + c_3 = 12$ implying $c_1 - c_3 \geq 9.$
This gives $c_1 \geq 9$ and so $c_2 \leq 3.$ Thus there is a vertex
$w \in V(G)$ with $S(w) \subset C_1.$ Since $\Delta(F) = 4$ exactly
three edges incident with $w$ are coloured with a colour from $C_3$
which implies $c_1 = 10, c_2 = c_3 = 1.$ So there are vertices $w_1,
w_2, w_3 \in V(G)$ such that $c(w_1w_2)=c(w_1w_3)=i \in C_2.$ Now
$S(w_2) \subset C_1$ and hence there are two edges incident with
$w_2$ which have a colour different from $S(w_2) \cup \{i\}.$ But
then $\Delta(F) \geq 5,$ a contradiction.

Let $w_1 \in V(F)$ be a vertex with $d_F(w_1) = 4$ and let $w_2,
\ldots, w_5$ be its neighbours. Let $w_6,w_7$ be the remaining
vertices. If $|N_F(w_i) \cap N_F(w_1)| \geq 3$ for a vertex $w_i$
with $i \in \{6,7\},$ then there is a rainbow $K_{2,3}$ as in the
previous case. Hence we may assume that $d_F(w_6) = d_F(w_7) = 3$
and $w_6w_7 \in E(F).$ Let $w_2w_6, w_3w_6 \in E(F).$ If  $w_2w_7,
w_3w_7 \in E(F),$ then there is a rainbow $K_{2,3}.$ Hence we may
assume that $w_2w_7 \notin E(F).$ With $|S(w_2)| \geq 3$ there is an
edge $w_2w_i$ for some $3 \leq i \leq 5$ with $c(w_2w_i)=p \in C^*.$
If $i \in \{4,5\},$ then $c(w_3w_i) \notin
\{c(w_1w_2),c(w_1w_3),c(w_2w_6),c(w_2w_7),c(w_2w_i),\}$ and so
$\{w_2, w_3, w_1, w_i, w_6\}$ form a rainbow $K_{2,3}.$ Otherwise we
have $N_F(w_7) \cap N_F(w_1) = \{w_4, w_5\}$ and $c(w_2w_3) = p,
c(w_4w_5) = q$ for some $p \neq q$ with $p,q \in C_1.$ Then
$c(w_iw_j) = r$ for some $r \in C_3$ for $2 \leq i \leq 3, 4 \leq j
\leq 5.$ Now observe that the four edges $w_6w_i, 2 \leq i \leq 5,$
are coloured with exactly two distinct colours. If $c(w_6w_4) =
c(w_6w_2),$ then $\{w_3,w_4,w_1,w_2,w_6\}$ form a rainbow $K_{2,3}.$
If $c(w_6w_4) = c(w_6w_3),$ then $\{w_2,w_4,w_1,w_3,w_6\}$ form a
rainbow $K_{2,3}.$ \qed \Qed

\noindent{$\bf rb(K_8, K_{2,3}) = 14$}\\
The proof is given in the appendix.

\section{Rainbow numbers for the house}

The complement $\overline{P_5}$ of the graph $P_5$ is known as the
{\it house} $H.$ Then

$${\cal H(H)} = \{C_5,C_4^+,B,Z_2\},$$

where $Z_2$ is the graph consisting of a $K_3$ with an attached path
with two edges.

Next observe that $C_3 \subset B, C_3 \subset Z_2, C_4 \subset C_4^+$
and therefore $ext(n, \{C_3, C_4, C_5\}) \leq ext(n, \{C_5, C_4^+,
B, Z_2\}).$ In \cite{GKL} projective planes have been used to
construct an infinite family of graphs with girth 6 (i.e. $\{C_3,
C_4, C_5\}$-free) having $\Omega(n^{\frac{3}{2}})$ edges. So we
deduce

\begin{cor}\label{c7}
$rb(K_n, H) = \Omega(n^{\frac{3}{2}}).$
\end{cor}

In order to show an upper bound we now show

\begin{thm}\label{t12}
$rb(K_n, H) \leq rb(K_n, K_{2,4})$ for all $n \geq 6.$
\end{thm}

\proof Let $K_n$ be edge coloured and let $F$ be a rainbow $K_{2,4}$
with vertices $w_1, w_2, \ldots, w_6$ and edges $w_iw_j$ for $1 \leq
i \leq 2, 3 \leq j \leq 6.$ Let $c(w_1w_{i+2})=i, c(w_2w_{i+2})=i+4$
for $1 \leq i \leq 4.$ Now consider the edge $w_3w_4.$ If $c(w_3w_4)
= x$ with $x \in C \setminus \{1,2, \ldots, 8\},$ then $F + e$
contains a rainbow $H.$ If $c(w_3w_4) = c(e) \in \{3,4,7,8\}$ for an
edge $e \in E(F),$ then $F + w_3w_4 - e$ contains a rainbow $H.$ If
$c(w_3w_4) = c(e) \in \{1,2,5,6\}$ for an edge $e \in E(F),$ then $F
+ w_3w_4 - e$ contains a rainbow $H.$

\qed

Corollary \ref{c7}, Theorem \ref{t11} and Theorem \ref{t12} yield

\begin{thm}
$rb(K_n, H) = \Theta(n^{\frac{3}{2}}).$
\end{thm}

\medskip

Now we will compute the rainbow numbers $rb(K_n, H)$ for $5 \leq n
\leq 8.$

\begin{lem}
$rb(K_{4+i},H) \geq 8 + {i \choose 2}$ for $1 \leq i \leq 4.$
\end{lem}

\proof Consider a $K_{4+i}$ and take a $K_4$ and a $K_i,$ which are
vertex disjoint. Colour the ${4 \choose 2} + {i \choose 2}$ included
edges distinct and all remaining edges with a new colour. Then there
is no rainbow $H$ and hence $rb(K_{4+i}) \geq ({4 \choose 2} + {i
\choose 2} + 1) + 1= 8 + {i \choose 2}.$ \qed

\noindent{$\bf rb(K_5, H) = 8$}\\
\proof Consider $8$ edges with distinct colours. The induced graph
is isomorphic to either $K_5 - P_3$ or $K_5 - 2K_2.$ In both cases
we find a rainbow $H.$ \qed

\noindent{$\bf rb(K_6, H) = 9$}\\
\proof Let the edges of $K_6$ be coloured with $9$ colours and let
$F$ be a rainbow subgraph of $K_6$ with $9$ edges. Since $rb(K_6,
C_5) = 9, F$ contains a rainbow $C_5.$ Let $w_1, w_2, \ldots, w_5$
be the vertices of the $C_5$ with edges $w_iw_{i+1}$ and colours
$c(w_iw_{i+1})=i$ for $1 \leq i \leq 5.$

If $c(e) \in \{6,7,8,9\}$ for a chord $e$ of the $C_5,$ then there
is a rainbow $H.$ Hence we may assume that $c(e) \in \{1,2, \ldots,
5\}$ for all chords $e$ of the $C_5.$ Let $w_6$ be the sixth vertex.
Then $d_F(w_6)=4.$ So we may assume that $w_6w_i \in E(F)$ with
$c(w_6w_i)=5+i$ for $1 \leq i \leq 4.$ Then we always find a rainbow
$H.$ \qed

\noindent{$\bf rb(K_7, H) = 11$}\\
\proof Let the edges of $K_7$ be coloured with $11$ colours and let
$F$ be a rainbow subgraph of $K_7$ with $11$ edges. We follow the
proof for $n = 6$ and consider the $C_5.$ Let $w_6$ and $w_7$ be the
two other vertices. If $d_{C_5}(w_i) \geq 4$ for $i = 6$ or $i = 7,$
then there is a rainbow house. Hence we may assume that $d_{C_5}
\leq 3$ for $i=6,7.$ If (i.e.) $d_{C_5}(w_6)=3$ and $w_6w_i,
w_6w_{i+1}, w_6w_{i+3} \in E(F)$ for some $i$ with $1 \leq i \leq
5,$ then there is a rainbow house. So we may assume that
$N_{C_5}(w_6) = \{w_i, w_{i+1}, w_{i+2}\} ,N_{C_5}(w_7) = \{w_j,
w_{j+1}, w_{j+2}\}$ for some $1 \leq i,j \leq 5,$ if
$d_{C_5}(w_6)=3$ or $d_{C_5}(w_7)=3.$ We distinguish two cases.

\medskip

\noindent{\bf Case 1} $d_{C_5}(w_6)=d_{C_5}(w_7)=3$\\
We have $N_{C_5}(w_6) = \{w_i, w_{i+1}, w_{i+2}\}$ for some $i$ with
$1 \leq i \leq 5.$ If $N_{C_5}(w_7) = \{w_i, w_{i+1}, w_{i+2}\}$ or
$N_{C_5}(w_7) = \{w_{i+1}, w_{i+2}, w_{i+3}\},$ then there is a
rainbow house. So we may assume that $N_{C_5}(w_7) = \{w_{i+2},
w_{i+3}, w_{i+4}\}.$ Then $c(w_6w_7)$ has the same (repeated) colour
as an edge from $F.$ Now observe that the graph $F + w_6w_7 - e$
contains a copy of $H$ for any edge $e \in E(F).$ Therefore, $G$
contains a rainbow house.

\medskip

\noindent{\bf Case 2} $d_{C_5}(w_6)=3, d_{C_5}(w_7)=2$\\
Then $w_6w_7 \in E(F).$ We have $N_{C_5}(w_6) = \{w_i, w_{i+1},
w_{i+2}\}$ for some $i$ with $1 \leq i \leq 5$ and $N_{C_5}(w_7) =
\{w_j, w_{j+1}\}$ or $N_{C_5}(w_7) = \{w_j, w_{j+2}\}$ for some $1
\leq j \leq 5.$ Now in all possible combinations we find a rainbow
house. \qed

\noindent{$\bf rb(K_8, H) = 14$}\\
\proof Let the edges of $K_8$ be coloured with $14$ colours and let
$F$ be a rainbow subgraph of $K_8$ with $14$ edges. If $|S(v)| \leq
3$ for a vertex $v \in V(G),$ then $G$ has a rainbow house by
induction. Hence we may assume that $|S(v)| \geq 4$ for every vertex
$v \in V(G).$ Consider the $C_5$ as above. Since five colours are
used for the edges of $G[V(C_5)],$ every vertex of this $C_5$ is
incident to at least two more edges in $F.$ Then $14 = |E(F)| \geq 5
+ 2 \cdot 5 = 15,$ a contradiction. \qed

\medskip

All computed rainbow numbers for the graphs $D, H, K_{2,3}$ are
listed in the following table.

\begin{center}
\begin{tabular}{c||c|c|c}
$n$ & $rb(K_n,D)$ & $rb(K_n, H)$ & $rb(K_n,K_{2,3})$\\
\hline\hline 5 & 7 & 8 & 8\\
6 & 8 & 9 & 10\\
7 & 10 & 11 & 12\\
8 & 12 & 14 & 14\\
9 & 14 && \\
10 & 17 &&
\end{tabular}
\end{center}

\noindent{\bf Acknowledgement:} We thank Maria Axenovich and Zsolt
Tuza for drawing our attention to rainbow numbers for bipartite
graphs and Jana Neupauerova for some stimulating discussions on the
computation of the rainbow numbers for the bull.

\section{Appendix}

\noindent{$\bf rb(K_8, K_{2,3}) = 14$}\\
\proof Consider two vertex disjoint $K_4's,$ colour all six edges of
one of them with colours $1, 2, \ldots, 6$ and all six edges of the
other one with colours $7, 8, \ldots, 12.$ All edges between the two
copies are coloured with colour $13.$ Then $G$ has no rainbow
$K_{2,3}.$ Hence $rb(K_8, K_{2,3}) \geq 14.$

Let the edges of $K_8$ be coloured with $14$ colours and let $F$ be
a rainbow subgraph of $K_8$ with $14$ edges. If $|S(v)| \leq 2$ for
a vertex $v \in V(G),$ then $G$ has a rainbow $K_{2,3}$ by
induction. Hence we may assume that $|S(v)| \geq 3$ for every vertex
$v \in V(G).$ Since $rb(K_8, D) = 12$ we may choose $F$ such that it
contains a rainbow $D.$ We distinguish two cases.

\noindent{\bf Case 1} $G$ contains a rainbow $K_4$\\
Let $F$ be a rainbow subgraph containing a $K_4$ with vertices $R =
\{w_1, w_2, w_3, w_4\}$ and let $T = \{w_5, w_6, w_7, w_8\}$ be the
other four vertices. If $|N_F(v) \cap R| \geq 2$ for a vertex $v \in
T,$ then there is a rainbow $K_{2,3}.$ Hence we may assume that
$|N_F(v) \cap R| \leq 1$ for every vertex $v \in T.$ With $|S(v)|
\geq 3$ we conclude that $d_T(w_i) \geq 2$ for $5 \leq i \leq 8.$
Hence $F[T]$ contains a $C_4$ implying that $F[T]$ is isomorphic to
$C_4$ or $D$ or $K_4.$ Now observe that in each of these three
possible cases we can find a cycle, say $w_5w_6w_7w_8w_5,$ such that
$|N_R(w_5)|=|N_R(w_7)|=1.$

If $w_iw_5, w_iw_7 \in E(F)$ for a vertex $w_i, 1 \leq i \leq 4,$
then there is a rainbow $K_{2,3}.$ Hence we may assume that $w_iw_5,
w_jw_7 \in E(F)$ for two vertices $w_i, w_j$ with $1 \leq i < j \leq
4,$ say $i=1, j=2.$ Consider the graph $F + w_1w_7.$ In this graph
$\{w_1, w_2, w_3, w_4, w_7\}$ and $\{w_5, w_7, w_1, w_6, w_8\}$
generate a $K_{2,3},$ where at least one of them is rainbow.

\medskip

\noindent{\bf Case 2} $G$ contains a rainbow $D$ but no rainbow $K_4$\\
Let $F$ be a rainbow subgraph containing a $D$ with vertices $R =
\{w_1, w_2, w_3, w_4\}$ and let $T = \{w_5, w_6, w_7, w_8\}$ be the
other four vertices. If $|N_F(w_i) \cap R| \geq 3$ for a vertex
$w_i, 5 \leq i \leq 8,$ then there is a rainbow $K_{2,3}.$ If
$|N_F(w_i) \cap R| = 2$ for a vertex $w_i, 5 \leq i \leq 8,$ then
there is a rainbow $K_{2,3}$ or a triangulated $C_5,$ which is
rainbow. We will use the abbreviation $TC_5$ for a triangulated
$C_5.$ Finally, if $|N_F(w_i) \cap R| = 1$ for $5 \leq i \leq 8$ and
$F[T] \cong D,$ then the four edges between $R$ and $T$ form a
perfect matching (otherwise we can find a rainbow $K_{2,3}$ or
$TC_5$ using previous arguments). We may assume that $w_iw_{i+4} \in
E(F)$ for $1 \leq i \leq 4.$ Consider the graph $F + w_1w_6.$ Then
$\{w_1, w_2, w_3, w_4, w_6\}$ and $\{w_1, w_5, w_6, w_7, w_8\}$
generate a subgraph with five vertices and seven edges, where at
least one of them is rainbow. This rainbow subgraph is either
isomorphic to $TC_5$ or contains a (rainbow) $K_{2,3}.$ Hence we may
assume that $F$ contains a rainbow $TC_5$ with vertices $R = \{w_1,
w_2, w_3, w_4, w_5\}$ and edges $\{w_1w_2, w_2w_3, w_3w_4, w_4w_5,
w_5w_1, w_1w_3, w_1w_4\}.$ Furthermore, let $c(w_1w_2)=1,
(w_2w_3)=2, (w_3w_4)=3, (w_4w_5)=4, (w_5w_1)=5, (w_1w_3)=6,
(w_1w_4)=7.$ Let $x = w_2w_5, y = w_2w_4, z = w_3w_5.$

If $c(x) \in \{1,5,6,7\},$ then there is a rainbow $K_{2,3}.$ So we
may assume $c(x) \in \{2,3,4\}.$

If $c(y) = 4,5$ then there is a rainbow $K_4^+,$ where $K_4^+$
denotes a $K_4$ with a pendant edge. If $c(y) = 2,7$ then there is a
rainbow $K_{2,3}.$ Hence we may assume that $c(y) \in \{1,3,6\}$ and
by symmetry $c(z) \in \{3,5,7\}.$ If $c(y)=6,c(z)=7,$ then there is
a rainbow $K_{2,3}.$ So we may assume $c(y) \neq 6$ or $c(z) \neq 7$
and we distinguish two cases.

\medskip
\noindent{\bf Subcase 2.1} $c(y)=6, c(z) \neq 7$ (the case $c(y)
\neq
6, c(z)=7$ is symmetric)\\
Since there is no rainbow $K_4$ we deduce $c(x)=4$ and so $c(z) \neq
3$ implying $c(z)=5.$

\noindent{\bf Subcase 2.2} $c(y) \neq 6, c(z) \neq 7$\\
If $c(y) = 1, c(z) = 5,$ then there is a rainbow $K_{2,3}.$ So
assume that $c(y) \neq 1$ implying $c(y) = 3.$ If $c(x) = 4$ then
there is a rainbow $K_{2,3}.$ If $c(x) = 2$ then there is a rainbow
$K_4^+.$ So we deduce that $c(x) = 3.$
Now if $c(z)=5$ then there is a rainbow $K_{2,3}.$ Hence we may
assume that $c(z)=3.$

Now for both subcases let $T = \{w_6, w_7, w_8\}$ be the three
remaining vertices. Since $|E(F)| = 14$ there is a vertex $w_i, 6
\leq i \leq 8,$ with $|N_R(w_i)| \geq 2.$ Let $u,w \in N_R(w_i).$ In
Subcase 2.1 every pair of vertices $w_i, w_j$ for $1 \leq i < j \leq
5$ is contained in a rainbow $4$-cycle with all its vertices from
$\{w_1, \ldots, w_5\}$ such that $w_i, w_j$ are not adjacent in this
cycle. This leads to a rainbow $K_{2,3}.$ In Subcase 2.2, if
$\{u,w\} \neq \{w_2,w_3\},\{w_4,w_5\},$ then there is a rainbow
$K_{2,3}.$ Hence we may assume that $w_2w_6, w_3w_6 \in E(F).$ Now
observe that the vertices $w_1$ and $w_3$ are both the center
(vertex of degree $4$) of a $TC_5.$ A repetition of the arguments
above gives $c(w_2w_4)=3$ and $c(w_2w_4)=1,$ a contradiction. \qed

\end{document}